\documentclass[12pt] {article}
\usepackage{amssymb,amsxtra,amsmath,amsthm}
\def \ZZ{{\mathbb{Z}}}

\def \CC{{\mathbb{C}}}
\def \FF{{\mathbb{F}}}
\def \PP{{\mathbb{P}}}

\def \OOO{{\mathcal{O}}}
\def \TTT{{\mathcal{T}}}

\def \aa{{\mathfrak a}}

\def \pp{{\mathfrak p}}
\def \qq{{\mathfrak q}}

\def\det{\mathop{\mathrm{det}}}

\def\mod{\mathop{\mathrm{mod}}}

\def\vert{\mathop{\mathrm{vert}}}
\def\edge{\mathop{\mathrm{edge}}}
\def\star{\mathop{\mbox{\Large$*$}}}

\def\ncs{\mathop{\mathrm{ncs}}}
\def\nncs{\mathop{\mathrm{nncs}}}
\def\gncs{\mathop{\mathrm{gncs}}}
\def\ngncs{\mathop{\mathrm{ngncs}}}

\def\m{\mathop{\mathrm{m}}}
\def\r{\mathop{\mathrm{r}}}

\def\ql{\mathop{\mathrm{ql}}}
\def\Gamma{\varGamma}
\def\Delta{\varDelta}

\begin{document}

\begin{center}
{\Large {\bf Genuine non-congruence subgroups of Drinfeld modular groups}}\\
\bigskip
{\tiny {\rm BY}}\\
\bigskip
{\sc A. W. Mason}\\
\bigskip
{\small {\it Department of Mathematics, University of Glasgow\\
Glasgow G12 8QW, Scotland, U.K.\\
e-mail: awm@maths.gla.ac.uk}}\\
\bigskip
{\tiny {\rm AND}}\\
\bigskip
{\sc Andreas Schweizer\footnote{The second author was supported
by the National Research Foundation of Korea (NRF) grant funded
by the Korean government (MSIP) (ASARC, NRF-2007-0056093).}}\\
\bigskip
{\small {\it Department of Mathematics,\\
Korea Advanced Institute of Science and Technology (KAIST),\\
Daejeon 305-701, South Korea \\
e-mail: schweizer@kaist.ac.kr}}\\
\end{center}
\begin{abstract}
\noindent
Let $A$ be the ring of elements in an algebraic function field $K$ over
a {\it finite} field $\FF_q$ which are integral outside a fixed place
$\infty$. In an earlier paper we have shown that the {\it Drinfeld
modular group} $G=GL_2(A)$ has automorphisms which map congruence
subgroups to non-congruence subgroups. Here we prove the existence
of (uncountably many) normal {\it genuine} non-congruence subgroups,
defined to be those which remain non-congruence under the action of
{\it every} automorphism of $G$. In addition, for all but finitely
many cases we evaluate $\ngncs(G)$, the smallest index of a normal
genuine non-congruence subgroup of $G$, and compare it to the minimal
index of an arbitrary normal non-congruence subgroup.
\\ \\
{\bf Key words:} Drinfeld modular group; genuine non-congruence subgroup;
non-standard automorphism
\\ \\
{\bf Mathematics Subject Classification (2010):}
11F06, 20E36, 20G30, 20H05
\end{abstract}

\subsection*{Introduction}

\noindent
An important point in the theory of the classical 
modular group is the free product decomposition
$$PSL_2(\ZZ)\cong \ZZ/2\ZZ\star\ZZ/3\ZZ.$$
As a consequence, every finite group $H$ that can be generated 
by an involution and an element of order $3$ can be obtained as 
quotient group $H\cong SL_2(\ZZ)/N$ for a normal finite index 
subgroup $N$ of $SL_2(\ZZ)$. Taking the compact Riemann surfaces 
that are the quotients of the extended complex upper half-plane
by $N$ resp. $SL_2(\ZZ)$ one obtains a Galois cover of $\PP^1(\CC)$
with prescribed Galois group $H$. Some of these quotients also 
furnish examples of Riemann surfaces whose automorphism group is 
big compared to their genus.
\par
Now for congruence subgroups $N$, the composition factors of 
$SL_2(\ZZ)/N$ can only be cyclic or $PSL_2(\FF_p)$. So for most
$H$ the group $N$ necessarily is a non-congruence subgroup.
\\ \\
For the analogue in positive characteristic, 
let $K$ be an algebraic function field of one variable with constant
field $k=\FF_q$. As usual we assume that $k$ is algebraically closed 
in $K$. Let $\infty$ be a fixed place of $K$ of degree $\delta$ and 
let $A$ be the ring of all those elements of $K$ which are integral 
outside $\infty$. The simplest example: If $K$ is a rational function
field over $k$ and $\delta=1$, then $A\cong k[t]$, a polynomial ring 
over $k$.
\par
Our focus of attention here are the {\it Drinfeld modular groups} 
$G=GL_2(A)$. These groups play a central role [G2] in the theory of 
{\it Drinfeld modular curves}, analogous to that of the classical 
modular group $SL_2(\ZZ)$ in the theory of modular curves and modular 
forms. 
\par
From the quotients of the
{\it Drinfeld upper half-plane} by finite index subgroups $N$ of $G$ one
obtains algebraic curves in positive characteristic. If $N$ is normal in
$G$ and contains the center $Z$ of $G$, then $G/N$ is a subgroup of the
automorphism group of that curve. Again, depending on what one would want
$G/N$ to be, in most cases $N$ has to be a non-congruence subgroup of $G$.
(See Lemma 2.5).
\par
In an earlier paper [MS3] we have shown that $SL_2(A)$ has
automorphisms which map (some) congruence subgroups to non-congruence
subgroups. Such automorphisms can be readily extended to $G$ (Theorem 2.2). 
They were originally introduced for the special case $A=k[t]$ by Reiner 
[R]. We extend a terminology used in [MS3].
\\ \\
\noindent {\bf Definition.}
We call an automorphism $\varPsi$ of an arithmetic group $X$
{\bf non-standard} if there exists a congruence subgroup $C$
of $X$ such that $\varPsi(C)$ is non-congruence.
\\
\noindent Otherwise $\varPsi$ is called {\bf standard}.
\\ \\
Non-standard automorphisms have been used in [MS3] to construct many 
non-congruence subgroups with prescribed properties. However, by 
construction, these non-congruence subgroups are subject to the same 
group-theoretic restrictions as congruence subgroups. In particular,
if such  a group $N$ is normal, the possible composition factors of
$G/N$ are known (Lemma 2.5). So if we want, just as an example, the 
simple group $A_7$ to be among these composition factors, then $N$
cannot be a congruence subgroup or obtained from one by a non-standard 
automorphism. What we need must be a ``genuine'' non-congruence subgroup
in the following sense.
\\ \\
\noindent {\bf Definition.}
Let $S$ be a non-congruence subgroup of an arithmetic group $X$. We call
$S$ {\bf genuine} if and only if $\varPsi(S)$ is a non-congruence subgroup,
for all $\varPsi \in \mathrm{Aut}(X)$.
\\ \\
Actually, using composition factors of $G/N$ that cannot occur for 
congruence subgroups, the existence of genuine non-congruence subgroups 
of Drinfeld modular group can be proved relatively easily from known 
results.
\\ \\
\noindent {\bf Theorem A.} \it
With every $G$ there corresponds a finite set of $A$-ideals 
$\mathcal{S}=\mathcal{S}(G)$ with the property that,
for all $\qq \notin \mathcal{S}$, there exist infinitely 
many normal genuine non-congruence subgroups $N$ for which
$$\ql(N)=\qq,$$
where $\ql(N)$ is the quasi-level of $N$.
\rm
\\ \\
The definition of the {\it quasi-level} of a subgroup of $G$ extends the 
classical definition of the {\it level} of a subgroup (originally applied 
to  the modular group $SL_2(\ZZ)$ or $PSL_2(\ZZ)$). The existence of 
uncountably many normal genuine non-congruence subgroups of any $G$ follows 
in similar way from an earlier result of Lubotzky [L, Theorem A(ii)] which 
states that $G$ always has a finite index subgroup mapping onto the free 
group of rank $2$. Lubotzky uses this result to prove [L, Theorem B(ii)] 
that $\hat{F}_{\omega}$, the free profinite group on countably many generators, 
is a closed subgroup of $C(SL_2(A))$, the {\it congruence kernel} of 
$SL_2(A)$.  The relationship between $\hat{F}_{\omega}$ and $C(SL_2(A))$ 
is further explored in [M3] and [MPSZ]. 
\\ \\
In previous papers [M4], [MS1], [MS2] we determined the smallest index
of a non-congruence subgroup of $SL_2(A)$. Here we turn our attention to
$\ngncs(G)$, the smallest index of a normal {\it genuine} non-congruence
subgroup of $G$. Our second main result (in Section 5) is the following.
\\ \\
\noindent {\bf Theorem B.} \it
In all but finitely many cases
$$\ngncs(G)=\m(G)=2,$$
\noindent where $\m(G)$ is the smallest index of a proper subgroup of $G$.
\rm
\\ \\
Note that this implies that the minimal index of a (not necessarily normal)
genuine non-congruence subgroup is then also $2$. We also determine $ngncs(G)$
for a number of cases for which Theorem B does not hold.
In particular, for the very important special case where $A$ is a polynomial
ring $\FF_q[t]$ we are able to show that $\ngncs(G)$ is strictly bigger than
the minimal index of a normal non-congruence subgroup of $G$.
See Section 6 for more and more precise results.
\\ \\
Although non-standard automorphisms exist for every Drinfeld
modular group this is not true for other important arithmetic groups like
$SL_2(\ZZ)$ and the Bianchi groups $SL_2(\OOO)$, where $\OOO$ is the ring
of integers in an imaginary quadratic number field. For these groups their
automorphisms are known [HR], [SV] to be ``standard''. So for them the 
definition of a genuine non-congruence subgroup is redundant, i.e. all
their non-congruence subgroups are genuine.
\par
Other known results in characteristic zero include the following. It is 
known [D, Proposition 4.3] that $SL_2(\ZZ)$ is a characteristic subgroup 
of $GL_2(\ZZ)$ which, in view of [HR], implies that every non-congruence 
subgroup of $GL_2(\ZZ)$ is genuine. It is worth noting that results for 
congruence subgroups of groups of type $SL_2$ and $GL_2$ do not necessarily 
apply to the corresponding projective groups $PSL_2$ and $PGL_2$. For example 
it is known that every non-congruence subgroup of $PSL_2(\ZZ)$ is genuine 
by [D, Corollary 4.4 (2)]. On the other hand, it is also known [JT] that 
$PGL_2(\ZZ)$ has a non-standard automorphism which therefore cannot fix 
$PSL_2(\ZZ)$. Hence $PSL_2(\ZZ)$ is not characteristic in $PGL_2(\ZZ)$.
\\

\subsection*{Notation}

\begin{tabular}{ll}
$k=\FF_q$ & the finite field of order $q=p^n$;\\
$K$         & an algebraic function field of one variable with
constant field $\FF_q$;\\
$g=g(K)$         & the genus of $K$;\\
$\infty$    & a chosen place of $K$;\\
$\delta$    & the degree of the place $\infty$;\\
$\nu$ & the discrete valuation of $K$ defined by $\infty$;\\
$\TTT$ & the Bruhat-Tits tree of $GL_2(K_\infty)$;\\
$A$ & the ring of all elements of $K$ that are integral outside $\infty$;\\
$G$ & the group $GL_2(A)$;\\
$\Gamma$ & the group $SL_2(A)$;\\
$Z$ & the centre of $G$;\\
$X$ & the groups $G,\Gamma$.\\
\end{tabular} \\

\subsection*{1. Subgroups defined by ideals}

Let $\qq$ be an $A$-ideal. It
is known that $A/\qq$ is {\em finite}, when $\qq$ is non-zero. We recall that by the well-known product formula $\nu(a) \leq 0$, for all $a \in A$,
and that $\nu(a)=0$ if and only if $ a \in k^*$.
\par

\noindent For each $\alpha, \beta \in k^*, a\in A$, we put
$$L(\alpha,\beta,a):=\left[\begin{array}{cc}\alpha&a\\
0&\beta\end{array}\right]$$
and
$$T(a):=L(1,1,a)= \left[\begin{array}{cc}1&a\\
0&1\end{array}\right].$$
\noindent For every subset $S$ of $A$ we put
$$T(S):=\{T(s):\; s\in S\}.$$
\noindent For each subgroup $H$ of $G$ and $g \in G$ we denote the conjugate
$gHg^{-1}$ by $H^g$. \\ \\
\noindent {\bf Definition.} We define the {\bf quasi-level} of $H$ to be
$$\ql(H):=\left\{h \in A:\; T(h) \in H^g,\; \mathrm{for\; all}\; g \in G\right\}.$$

\noindent The {\bf level} of $H$ is the biggest $A$-ideal contained in $\ql(H)$.\\ \\
\noindent Clearly $\ql(H)$ is an additive subgroup of $A$ and, by considering
conjugates by the elements $\mathrm{diag}(\alpha,1)$, where $\alpha \in k^*$,
it is clear that it is also a vector space over $k$.\\ \\
In [MS3] the ``quasi-level" of $H\leq\Gamma$ is defined to be
$$\ql(H)^*:=\left\{h \in A:\; T(h) \in H^g,\; \mathrm{for \;all}\; g \in \Gamma\right\}.$$
\noindent It is easily shown that the definitions are equivalent.
\\ \\
\noindent {\bf Lemma 1.1.} \it
$\ql(H)^*=\ql(H).$
\rm
\\ \\
\noindent {\bf Proof.}
Clearly $\ql(H) \leq \ql(H)^*$. By conjugating with the matrices $\mathrm{diag}(\alpha,\alpha^{-1})$ and using the fact that every element in a finite field is a sum of two squares,
it follows that $\ql(H)^*$ is a vector space over $k$.\\
\noindent Now let $h \in \ql(H)^*$ and $g \in G$. Then $g=dg'$, where $d =\mathrm{diag}(\beta,1)$ (with $\beta \in k^*$) and $g' \in \Gamma$. Then, from the above, $T(h \beta^{-1}) \in H^{g'}$ and so $T(h) =dT(h \beta^{-1})d^{-1} \in H^g$.
\hfill $\Box$
\\ \\
\noindent We define the Borel subgroup of $G$
$$B_2(A)=G_{\infty}=\left\{ L(\alpha,\beta,a):\;\alpha,\beta \in k^*,a \in A\right\}.$$
\par
\noindent {\bf Definition.} For every ideal $\qq$ of $A$ we define
$$\Delta(\qq):= \mathrm{ the\; normal\; subgroup\; of\;} \Gamma \;\mathrm
{generated\; by}\; T(\qq).$$
\noindent A subgroup $H$ of $G$ is said to have {\it non-zero level} if
$\Delta(\qq)\subseteq H$, for some $\qq\neq\{0\}$. Otherwise $H$ is said
to have {\it level zero}.
\\ \\
As in the proof of Lemma 1.1 one shows that $\Delta(\qq)$ is also normal
in $G$. So $\Delta(\qq)$ can also be defined as the normal subgroup of $G$
generated by $T(\qq)$.
\\ \\
\noindent {\bf Definition.} Let
$$\Gamma(\qq)=\left\{M \in \Gamma: M \equiv I_2\;(\mod\;\qq)\right\}.$$
\noindent A subgroup $C$ of $G$ is said to be a {\bf congruence subgroup}
if $\Gamma(\qq') \leq C$, for some $\qq' \neq \{0\}$. Such a subgroup is
necessarily of finite index in $G$. A finite index subgroup of $G$ which
is not congruence is a called a {\bf non-congruence subgroup}.
\\ \\
\noindent We will make use of the following properties of these subgroups. The first
[B, (9.3) Corollary, p.267] plays an important role in determining whether
or not a finite index subgroup of $G$ is congruence (the so-called
{\it congruence subgroup problem}).
\\ \\
\noindent {\bf Lemma 1.2.} \it
Let $\qq_1,\qq_2$ be $A$-ideals, where $\qq_2 \neq\{0\}$. Then
$$\Delta(\qq_1).\Gamma(\qq_2)=\Gamma(\qq_1+\qq_2).$$
\rm
\\
From this one easily obtains the following classical criterion 
for being a congruence subgroup.
\\ \\
{\bf Corollary 1.3.} \it
Let $H$ be a subgroup of $G$ with $\Delta(\qq)\subseteq H$ 
for some non-zero ideal $\qq$ of $A$. Then $H$ is a congruence 
subgroup if and only if $H$ contains $\Gamma(\qq)$.
\rm
\\ \\
{\bf Proof.}
If $H$ is a congruence subgroup, then by definition it contains
$\Gamma(\qq')$ for some non-zero ideal $\qq'$. So by Lemma 1.2 it 
contains $\Gamma(\qq+\qq')$, which contains $\Gamma(\qq)$.
The converse is clear.
\hfill $\Box$
\\ \\
In view of this criterion, the existence of non-congruence subgroups
of non-zero level is a consequence of the following result.
\\ \\
\noindent  {\bf Lemma 1.4.} \it
There exists an epimorphism
$$\Gamma(\qq)/\Delta(\qq)\twoheadrightarrow F_r,$$
\noindent where $F_r$ is the free group of (finite) rank $r=\r(\qq)=\mathrm{rk}_{\ZZ}(\Gamma(\qq))$, the torsion-free rank of the abelianization of $\Gamma(\qq)$. Moreover
$$\r(\qq)\rightarrow \infty\;\;as\;\; A/\qq \rightarrow \infty.$$
\rm
\\
\noindent {\bf Proof.}
Let $\varLambda(\qq)$ be the subgroup of $G$ generated by $\Gamma(\qq)\cap G_v$, for all $v \in \vert(\TTT)$. Then $\varLambda(\qq)\unlhd G$ and, from the theory of groups acting on trees [Se, Corollary 1, p.55], it is known that
$$\Gamma(\qq)/ \varLambda(\qq) \cong F_{\r(\qq)},$$
\noindent the fundamental group of the quotient graph $\Gamma(\qq) \backslash \TTT$, which is known [Se, Corollary 4, p.108]
to have {\it finite} rank $\r(\qq)$. Again from the presentation [Se, p.42] of $\Gamma(\qq)$ derived from its action on $\TTT$, together with [Se, Proposition 2, p.76], it follows that $\varLambda(\qq)$ is the subgroup of $\Gamma(\qq)$ generated by its torsion elements. Hence $\Delta(\qq) \leq \varLambda(\qq)$. \\
\noindent Given $\qq_1$ and any $\qq_2$ for which $\qq_2 \leq \qq_1$, it is clear that $F_{\r(\qq_2)} \leq F_{\r(\qq_1)}\cap\Gamma(\qq_2)$. By means of the Schreier formula we can always choose $\qq_2$ such that $\r(\qq_2)>\r(\qq_1)$.
\hfill $\Box$
\\

\subsection*{2. Non-standard automorphisms}

\noindent For each non-negative integer $n$, we put
$$A(n):=\{x\in A:\quad \nu(x)\geq -n\}.$$
\noindent Then $A(n)$ is a {\em finite-dimensional} vector space over
$k$, whose dimension is determined by the {\it Riemann-Roch
Theorem}. (See [St, Theorem 1.5.15, p.30].) Obviously $k\subseteq A(n)$.
In particular $A(0)=k$, by [St, Corollary 1.1.20, p.8].\\
\noindent We put
$$G_n:=\{L(\alpha,\beta,a):\alpha,\beta\in k^*,\; a\in A(n)\}.$$
\noindent Then $G_n$ is a finite subgroup of $G$. We note that
$$ G_{\infty}=\displaystyle{\bigcup_{n \geq 0}\;G_n}.$$
\\
\noindent Serre's decomposition theorem [Se, Theorem 10, p.119] shows that
$G_{\infty}$ is a (non-trivial) factor in a decomposition of
$G$ as an amalgamated product of a pair of its subgroups. See [MS3, Theorems 2.1, 2.2].
\\ \\
\noindent {\bf Theorem 2.1. } \it
Let $n$ be the smallest non-negative integer for which
$$\delta n\geq 2g-1.$$
\noindent Then there exists a subgroup $H$ of $G$, such that
$$G=G_\infty \star_{\quad{L}} H,$$
where $L=G_{n_0}$. Moreover
$$\mathrm{dim}_k(A(n_0))=n_0\:\delta+1-g.$$
\rm
\\
Note that $n_0=0$ when $g=0$. Otherwise $n_0 >0$.
\\ \\
\noindent {\bf Definition.} Let $\phi: A \rightarrow A$ by any
$k$-automorphism of the $k$-vector space $A$ which fixes the elements
of $A(n_0)$ (including $k$).  Then $\phi$ induces an automorphism
$\varPhi$ of $G_{\infty}$ defined by
$$\varPhi:L(\alpha,\beta,a)\longmapsto L(\alpha,\beta,\phi(a)),$$
where $\alpha,\beta\in k^*, a \in A$.
\\ \\
\noindent The following is an immediate consequence of Theorem 2.1.
\\ \\
{\bf Theorem 2.2. } \it
Let $\phi$ be any $k$-automorphism of $A$ which fixes the
elements of $A(n_0)$. The map
$$
\varPhi(g)=\left\{\begin{array}{cll}
L(\alpha,\beta,\phi(a))&, & g=L(\alpha,\beta,a)\in G_{\infty}\\[8pt]
g&, & g\in H
\end{array}\right.
$$
extends to an automorphism of $G$ (and $\Gamma$).
\rm
\\ \\
\noindent As shown in [R] and [MS3] in most cases such an automorphism is
{\it non-standard}. {\it Standard} automorphisms include inner automorphisms, 
the contragredient map $M \mapsto (M^T)^{-1}$, twists with certain determinant
characters, i.e. $M\mapsto \chi(\det(M))M$ where $\chi:k^*\to k^*$ is a group 
homomorphism (with the property that $\chi(\alpha^2)=\alpha^{-1}$ if and only 
if $\alpha=1$), or automorphisms derived from {\it ring}-automorphisms of 
$A$. Such a ring-automorphism $\psi$ of $A$ induces the automorphism 
$\varPsi$ of $G$, defined by
$$
\varPsi:\left[\begin{array}{ll}a&b\\
c&d\end{array}\right]\mapsto\left[\begin{array}{ll}\psi(a)&\psi(b)\\
\psi(c)&\psi(d)\end{array}\right].
$$
\noindent Clearly every such ring-automorphism maps $A$-ideals to $A$-ideals.
\par
Also note that the standard automorphism of $GL_2(\FF_q[t])$ induced by
the ring automorphism $t\mapsto t+1$ of $\FF_q[t]$ can also be obtained
via the construction in Theorem 2.2 from a suitable $\FF_q$-vector space 
automorphism $\phi$. This shows that not every automorphism in 
Theorem 2.2 is necessarily non-standard.
\par
For the case $A=\FF_q[t]$ it is known [R] that the standard automorphisms 
listed above together with the automorphisms from Theorem 2.2 generate 
$\mathrm{Aut}(G)$. This is the only case for which $\mathrm{Aut}(G)$ is 
known.
\par
We record two obvious properties.
\\ \\
\noindent {\bf Lemma 2.3.} \it
For every subgroup $S$ of $G$ and any automorphism 
$\varPhi$ of $G$ defined by $\phi$ (as in Theorem 2.2) we have
$$ \ql(\varPhi(S))=\phi(\ql(S)).$$
\noindent In particular
$$\ql(\varPhi(S))=A\Leftrightarrow \ql(S)=A.$$
\rm
\\
Note however that an {\it arbitrary} automorphism of $G$ or $\Gamma$ need 
not map subgroups of quasi-level $A$ to subgroups of quasi-level $A$.
\\ \\
\noindent {\bf Example 2.4.} Let
$$A=\FF_2[x,y]\ \ \ \hbox{\rm with}\ \ \ y^2 +xy=x^3 +x^2 +x.$$
Then by the results of Takahashi [T] we have
$$GL_2(A)\cong\Delta(\infty)\star\Delta(0)\star\Delta(1)$$
in the notation of [MS2, Theorem 5.3] (not of Section $1$ of the 
current paper). Explicitly, $\Delta(1)\cong\ZZ/3\ZZ$ and
$$\Delta(\infty)=GL_2(\FF_2)\star_{B_2(\FF_2)} B_2(A)$$
and $\Delta(0)$ is isomorphic (as an abstract group) to $\Delta(\infty)$.
So there is a an automorphism $\sigma$ of $GL_2(A)$ that interchanges
$\Delta(0)$ and $\Delta(\infty)$. 
\par
From the free product we can construct normal finite index subgroups $H$ 
containing $\Delta(0)$ and $\Delta(1)$ with practically any prescribed 
quasi-level. Then $\sigma(H)$ contains $\Delta(\infty)$, and hence has
quasi-level $A$.
\\ \\
Before proceeding we require a further definition.
\\ \\
\noindent {\bf Definition.} For each ideal $\qq$ let
$$Z(\qq)=\left\{X \in G: X\equiv \alpha I_2\;(\mod\;\qq)\;for \; some\;\alpha \in A\right\}.$$

\noindent It is clear that $Z(\qq)/\Gamma(\qq^2)$ is abelian.
\\ \\
\noindent {\bf Lemma 2.5.} \it
Let $N$ be a normal congruence subgroup of index $n$ in $G$ and let $\varPsi$ be any automorphism of $G$.
Then the (simple) factors in a composition series of $G/\varPsi(N)$ are either cyclic of prime order or
are isomorphic to some $PSL_2(\FF_{q^s})$, where $s \geq 1$.
\par
Moreover, if the level of $N$ is divisible by a prime $\pp$ with $|A/\pp|>3$,
then at least one factor in this composition series is of the latter type.
\rm
\\ \\
\noindent {\bf Proof.}
It is known [M2, Theorem 3.14] that, for some non-zero ideal $\qq_0$
 $$\Gamma(\qq_0^2\aa^2) \leq N \leq Z(\qq_0),$$
 \noindent where $\aa$ is the product of all prime ideals of $A$ of index $2$ or $3$. If none exists we put $\aa=A$. We may confine our attention therefore to the composition factors of the groups $Z(\qq)/\Gamma(\qq^2\aa^2)$ and $Z(\qq)/Z(\qq\pp)$, where $\pp$ as usual is prime.
\\
\noindent We write $\aa=\aa_1\aa_2$, where $\aa_1+\qq=A$ and $\aa_2$ divides $\qq$. From standard results (like Lemma 1.2) it follows that
$$\Gamma(\qq)/\Gamma(\qq^2\aa^2) \cong (\Gamma(\qq)/\Gamma(\qq^2\aa_2^2))\times \displaystyle{\prod_{\pp|\aa_1}\Gamma/\Gamma(\pp^2)}.$$
\noindent Now the first group in the decomposition is metabelian from above and each $\Gamma/\Gamma(\pp^2)$ is a soluble group of order $m^4(m^2-1)$, where $m=2,3$.
\\
\noindent We now consider the group $\Gamma(\qq)/\Gamma(\qq\pp)$. If $\pp$ divides $\qq$ it is abelian. On the other hand if $\pp+\qq=A$ then
$$\Gamma(\qq)/\Gamma(\qq\pp) \cong \Gamma/\Gamma(\pp) \cong SL_2(\FF_{q^s}),$$
\noindent where $\FF_{q^s}=A/\pp$. The first part follows. 
\par
Under the condition of the second statement we have 
$$N \leq Z(\pp)$$
again by [M2, Theorem 3.14]. We note that 
$$G/Z(\pp) \hookrightarrow PGL_2(\FF_{q^s}),$$
\noindent where $\FF_{q^s}=A/\pp$. It is well-known that $PSL_2(\FF_{q^s})$ is contained in this embedding. When $q^s >3$ the latter group is simple. The second part follows.
\hfill $\Box$
\\ \\
\noindent If $q>3$, the condition in the second statement of Lemma 2.5 is
of course automatic for all congruence subgroups except those containing 
$\Gamma$. We record a restricted version of this lemma.
\\ \\
{\bf Lemma 2.6.} \it
Let $N$ be a proper normal congruence subgroup of $\Gamma$ and let $\varPsi$ be any automorphism of $\Gamma$.
Then the factors in a composition series of $\Gamma /\varPsi(N)$ are as in Lemma 2.5.
\noindent Moreover if the level of $N$ is divisible by a prime $\pp$ with
$|A/\pp|>3$, then at least one factor in this composition series is of the 
latter type.
\rm
\\

\subsection*{3. Genuine non-congruence subgroups}

\noindent {\bf Notation.} For the remainder of this paper $X$ will always denote $G$ or $\Gamma$.
\\ \\
\noindent {\bf Definition.}  A finite index subgroup $S$ of $X$ is said
to be a {\bf genuine non-congruence subgroup} of $X$ if $\varPsi(S)$ is
a non-congruence subgroup, for all $\varPsi \in \mathrm{Aut}(X)$.
\\ \\
\noindent The following straightforward result enables us in many instances to assume that a given genuine non-congruence subgroup is normal.
\\ \\
\noindent {\bf Lemma 3.1.} \it
A finite index subgroup of $X$ is a genuine non-congruence subgroup if
and only if its core in $X$ is a genuine non-congruence subgroup of $X$.
\rm
\\ \\
\noindent {\bf Remark 3.2.}
Note however that we cannot be sure whether a genuine non-congruence
subgroup $H$ of $\Gamma$ automatically is a genuine non-congruence subgroup
of $G$, as theoretically $G$ might have non-standard automorphisms (other
than those discussed in Theorem 2.2) that do not respect $\Gamma$.
\\
\noindent
Similarly, if $N$ is a (normal) genuine non-congruence subgroup of $G$,
we cannot be sure whether $N\cap\Gamma$ is a genuine non-congruence
subgroup of $\Gamma$, as theoretically there might be automorphisms of
$\Gamma$ that do not extend to automorphisms of $G$.
\\ \\
\noindent We now make use of Lemma 1.4 to prove the existence of genuine
non-congruence subgroups.
\\ \\
\noindent {\bf Theorem 3.3.} \it
For all but finitely many $\qq$ there exist infinitely many
normal genuine non-congruence subgroups $N$ of $X$ for which
$$\ql(N)=\qq.$$
\rm
\\
\noindent {\bf Proof.} By Lemma 1.4 there exists an epimorphism
$$ \Gamma(\qq)/\Delta(\qq)\twoheadrightarrow F_2,$$
for all but finitely many $\qq$. We note that if $M$ is a finite
index subgroup of $X$ and
$$\Delta(\qq) \leq M \nleqq \Gamma(\qq),$$
\noindent then $M$ is non-congruence by Corollary 1.3.\\
\noindent Let $H$ be any finite group that can be generated by 
$2$ elements and that has a non-cyclic simple composition factor 
that is not isomorphic to any $PSL_2(\FF_{q^s})$. For example, we 
can always take $H=S_n$ with $n\geq 7$.
Or if $q$ is different from $2$, $4$, $5$, we can even take 
$H=A_5\cong PSL_2(\FF_5)\cong PSL_2(\FF_4)$. 
Then there exists $M \leq \Gamma(\qq)$, where $M \geq \Delta(\qq)$, 
for which
$$\Gamma(\qq)/ M \cong H.$$
\noindent Let $N$ be the core of $M$ in $X$. Then
$\Delta(\qq) \leq N \leq \Gamma(\qq)$ so that $\ql(N)=\qq$.
In addition $N$ is genuine by Lemmas 2.5, 2.6.
\hfill $\Box$
\\ \\
\noindent The subgroups in Theorem 3.3 have all non-zero level.
By an earlier method we can prove the following.
\\ \\
\noindent {\bf Theorem 3.4.} \it
There exist uncountably many normal genuine non-congruence subgroups
of $X$ of level zero.
\rm
\\ \\
\noindent {\bf Proof.} As in the proof of Theorem 3.3 we choose $\qq$ and
$N \unlhd X$ so that $\Delta(\qq) \leq N \leq \Gamma(\qq)$ and $A_7$, say,
is a factor in the composition series of $X/N$. Now choose an ideal $\qq_0$
so that $A=V_0\oplus \qq_0$, where $V_0$ is a finite-dimensional space
containing $A(n_0)$ (from Theorem 2.1). \\
\noindent Now, if $\qq_1 \leq \qq_2$ and $\qq_1 \neq \{0\}$, then the natural map
$$\Gamma(\qq_1) /\Delta(\qq_1)\longrightarrow \Gamma(\qq_2)/\Delta(\qq_2)$$
\noindent is {\it surjective} by Lemma 1.2. So replacing $\qq$ with
$\qq\cap\qq_0$ we may assume that $\qq=\qq_0$. Let $W$ be one of the
uncountably many subspaces of $A$ not containing any non-zero $A$-ideal
for which $A=V_0 \oplus W$. Then we can find a non-standard isomorphism
$\varPhi$ of $X$ for which $\ql(\varPhi(N))=W$.
\hfill $\Box$
\\ \\
\noindent  However as we now show many genuine non-congruence subgroups of $X$ of quasi-level $A$ do exist. Let
$$X_V=\left<X_v:\;v \in \vert(\TTT)\right>.$$
\noindent By [Se, Proposition 2, p.76] $X_V$ is the subgroup of $X$ generated by all its torsion elements and so is invariant under every automorphism of $X$. It follows that $\Delta(A) \leq X_V$. In addition
$$X/X_V \cong F_{\r(X)},$$
\noindent where $F_{\r(X)}$ is the fundamental group of the quotient graph $X \backslash \TTT$.
See [Se, Corollary 1, p.55]. The rank $\r(X)$ is known to be finite [Se, Corollary 4, p.108]. In particular $\r(X)=0$ if and only if $X\backslash \TTT$ is a {\it tree}. Moreover
 there are formulae for $\r(X)$ involving $\delta,q$ and values of the
{\it L-polynomial} of $K$ ([G1], [G2, p.73], or see [MS1, p.33].) From these
it is clear that, for any fixed
$g,\;  \r(X)\rightarrow \infty$, as $\delta,q\rightarrow\infty$.\\
\noindent It is clear that $\r(\Gamma) \geq \r(G)$. The rank zero cases are known
precisely [MS1, Theorem 2.10]. For convenience we record them.
\\ \\
\noindent {\bf Theorem 3.5.} \it
\begin{itemize}
\item[(i)] $\r(G)=0$ if and only if $(g,\delta)=(1,1),(0,1),(0,2)$ or $(0,3)$.
\item[(ii)] $\r(\Gamma)=0$ if and only if $(g,\delta)=(0,1), (0,2)$ or
(when $q$ is even) $(0,3),(1,1)$.
\end{itemize}
\rm
\bigskip
\noindent {\bf Lemma 3.6.} \it
Let $S$ be any proper finite index subgroup of $X$ containing $X_V$.
Then $S$ is a genuine non-congruence subgroup of $X$ with $\ql(S)=A$.
\rm
\\ \\
\noindent {\bf Proof.}
Recall that $\varPsi(X_V)=X_V$ for all $\varPsi \in \mathrm{Aut}(X)$.
Any congruence subgroup containing $X_V$ must contain $\Gamma.X_V=X$
by Lemma 1.2.
\hfill $\Box$
\\ \\
\noindent {\bf Corollary 3.7.} \it
Suppose that $r=\r(X)>0$ and that $H$ is any finite group with at most
$r$ generators. Then there exists a normal genuine non-congruence subgroup
$N$ of $X$ of level $A$ with
$$X/N \cong H.$$
\rm

\noindent When $\r(X)=0$ there need not be any non-congruence subgroups of
level $A$. Consider, for example, the case where $(g,\delta)=(0,1)$. Then
$A =k[t]$, a euclidean ring. In this case then $\Delta(A)=\Gamma$ and
$X_V=X$.
\\ \\
\noindent  On the one hand, Lemmas 2.5, 2.6 provide necessary conditions
for a non-congruence subgroup to be not genuine. On the other hand,
Corollary 3.7 enables us to show that these conditions are not sufficient.
\\ \\
\noindent {\bf Example 3.8.}
We provide a simple illustration of Corollary 3.7.
Consider one of the many $G$ with $\r(G) \geq 4$.
Let $\pp$ be any prime $A$-ideal. Then, by Lemma 1.2,
$$G/ \Gamma(\pp)\cong SL_2(A/\pp)\rtimes k^*.$$
\noindent Now $A/\pp$ is a field $\FF_{q'}$, for some power $q'$ of
$q$, and so $SL_2(A/\pp)$ is generated by all $T(a)$ and $T(a)^T$.
Let $\lambda$ be a generator of $\FF_{q'}^*$. Then $G / \Gamma(\pp)$
is generated by
$T(1)$, $\mathrm{diag}(\lambda,\lambda^{-1})$, $\mathrm{diag}(\mu,1)$,
where $\mu$ generates $k^*$ and
$\left[\begin{array}{cc}0&-1\\ 1&0\end{array}\right]$.
(Recall that every element in a finite field is a sum of $2$ squares.)
Then there exists $N\trianglelefteq G$, with $G_V \leq N$, such that
$$G/N \cong G/\Gamma(\pp).$$
\noindent In some cases, of course, the rank restriction here can be weakened.
For example, if $q'=q$, the group $G/ \Gamma(\pp)$ is $3$-generated.
And if $A/\pp \cong\FF_p$, a prime field, then 
$\Gamma/\Gamma(\pp)\cong SL_2(\FF_p)$ is even $2$-generated 
(as a quotient of $SL_2(\ZZ)$).
\\

\subsection*{4. Some immediate criteria for being genuine}

\noindent In this section we show that sometimes the index of a subgroup in itself can show that it is genuine non-congruence.
\\ \\
\noindent {\bf Proposition 4.1.}  \it
Let $N$ be a proper normal subgroup of index $n$ in $G$, where
$\mathrm{gcd}(n,q)=1$. Suppose further that there exists $S \leq N$ such that
$$S \cong k^* \times k^*.$$
\noindent Then $N$ is a genuine non-congruence subgroup of $G$.
\rm
\\ \\
\noindent {\bf Proof.}
Clearly $\Delta(A)\leq\varPsi(N)$ for any automorphism $\varPsi$ of $G$.
By Maschke's theorem applied to (the abelian group) $\varPsi(S)$, there
exists $ g \in GL_2(K)$ such that
$$g\varPsi(S)g^{-1}=D,$$
\noindent where $D$ is the set of all diagonal matrices in $GL_2(k)$.
Hence $\varPsi(N).\Gamma=G$ and so $\varPsi(N)$ is a non-congruence subgroup
of $G$ by Corollary 1.3.
\hfill $\Box$
\\ \\
\noindent {\bf Remark 4.2.}
If $q$ is even or if $n>2$, the condition $S\leq N$ in Proposition
4.1 can be replaced by $Z\leq N$ as then $Z\leq\varPsi(N)$ for any
$\varPsi\in\mathrm{Aut}(G)$. However, it is not clear whether conditions
like $\det(N)=k^*$ or $N.\Gamma=G$ would suffice, as their behaviour
under $\varPsi$ is not obvious.
\\ \\
\noindent The version of Proposition 4.1 for $\Gamma$ is simpler.
\\ \\
\noindent {\bf Proposition 4.3.}  \it
Let $N$ be a proper normal subgroup of index $n$ in $\Gamma$, where
$\mathrm{gcd}(n,q)=1$. Then $N$ is a genuine non-congruence subgroup
of $\Gamma$.
\rm
\\ \\
The following two results are easy conclusions of Lemmas 2.5 and 2.6.
\\ \\
\noindent {\bf Proposition 4.4.} \it
Suppose that $q>3$ and that $N$ is a normal subgroup of index $n$ in $G$,
where $n \nmid(q-1)$. If $|PSL_2(\FF_q)| \nmid n$, then $N$ is a genuine
non-congruence subgroup of $G$.
\rm
\\ \\
\noindent {\bf Proposition 4.5.} \it
Suppose that $q>3$ and that $N$ is a proper normal subgroup of index
$n$ in $\Gamma$. If $|PSL_2(\FF_q)| \nmid n$, then $N$ is a genuine
non-congruence subgroup of $\Gamma$.
\rm
\\ \\
\noindent Note that Propositions 4.4, 4.5 hold in particular if
$\mathrm{gcd}(n,q)=1$ or $\mathrm{gcd}(n,q\pm 1)=1$.\\
\noindent The restrictions
on $q$ are necessary. When $q=2,3$ and $A=k[t]$ it is well-known that
$X$ has normal congruence subgroups of index $q$. Moreover it is known
[MS1, Lemma 3.1] that for these cases $X$ has non-congruence subgroups
of index $q$ which are {\it not} genuine.
\\ \\
{\bf Notation.}
For the case where a group $H$ has proper finite index subgroups, we
denote by $\mathrm{m}(H)$ ($>1$) the smallest index of such a subgroup.
\\ \\
\noindent It is a classical result (originally due to Galois) that
$\mathrm{m}(SL_2(\FF_q))=q+1$ for $q>11$ and $q=4,8$. Otherwise this
index is $q$ unless $q=9$ in which case it is $6$.
\\ \\
\noindent For non-normal subgroups of $\Gamma$ we can now prove the following.
\\ \\
\noindent {\bf Proposition 4.6.} \it
Suppose that $q>3$ and that $H$ is a proper subgroup of $\Gamma$ for which
$$|\Gamma:H| < m(SL_2(\FF_q)).$$
\noindent Then $H$ is a genuine non-congruence subgroup of $\Gamma$.
\rm
\\ \\
\noindent {\bf Proof.}
Let $S=SL_2(\FF_q)$. Then, for each $g \in \Gamma$,
$$|S:S\cap H^g| \leq |\Gamma:H^g|.$$
\noindent It follows that $S \leq H^g$, and hence that $S$ is contained
in the core of $H$. By [MS3, Lemma 3.2] this implies $\varDelta(A) \leq H$.
By Corollary 1.3, $H$ then is non-congruence. We can repeat the argument with
$\varPsi(H)$, for all $\varPsi \in \mathrm{Aut}(\Gamma)$.
\hfill $\Box$
\\

\subsection*{5. The minimum index of a genuine non-congruence subgroup}

\noindent The first two of the following appear in [MS1] and [MS2].
\\ \\
\noindent {\bf Definitions.}
\begin{itemize}
\item[(i)] $ \ncs(X)= \mathrm{min}\{|X:S|: \; S \leq X, \;S\; noncongruence \}.$
\item[(ii)] $ \nncs(X)= \mathrm{min}\{|X:S|: \; S \leq X, \;S\; normal, \; noncongruence \}.$
\item[(iii)] $\gncs(X)=\mathrm{min}\{|X:S|:\; S \leq X,\;S\;genuine\; noncongruence \}.$
\item[(iv)] $\ngncs(X)=\mathrm{min}\{|X:S|:\; S \leq X,\;S\;normal,\; genuine\; noncongruence \}.$
\end{itemize}
\noindent In [M4], [MS1], [MS2] we determined $\ncs(\Gamma)$ in all cases,
and also $\nncs(\Gamma)$ [MS2, Theorem 6.2]. In this section we evaluate
$\ngncs(X)$ and $\gncs(X)$ in all but finitely many cases. It is clear that
$\ngncs(X)\geq\nncs(X)$ and that $\gncs(X)\geq\ncs(X)$.
\\ \\
An immediate consequence of Corollary 3.7 is the following.
\\ \\
\noindent {\bf Theorem 5.1.} \it
Suppose that $\r(X)>0$. Then
$$\ngncs(X)=\gncs(X)=\nncs(X)=\ncs(X)=\m(X)=2.$$
\rm

\noindent Theorem 5.1 also holds for some but not all rank zero cases,
which are listed in Theorem 3.5.  For the remainder of this section we consider the cases $(g,\delta)=(1,1),\;(0,3)$ in detail. We recall from [MS4] the possible structures
of the stabilizers in {\it any} $G$ of the vertices of $\TTT$. For any $v \in \vert(\TTT)$ it is known that one of the following holds
\begin{itemize}
\item[(a)] $G_v \cong GL_2(k)\; \mathrm{or}\; \FF_{q^2}^*$.
\item[(b)] $G_v \cong k^* \times N$,
\item[(c)] $G_v/N \cong k^* \times k^*$,
\end{itemize}
\noindent where $N \cong V^+$, the additive group of a finite dimensional $k$-space $V$. See [MS4, Corollaries 2.2, 2.4, 2.7]. We require the following.
\\ \\
\noindent {\bf Lemma 5.2.} \it Suppose that $(g,\delta)=(1,1)$ or $(0,3)$.
Then there exist subgroups $P,\;Q$ of $G$ such that
$$G=P\star_{\quad Z}Q,$$
\noindent where
 \begin{itemize}
 \item[{(i)}]$GL_2(k) \leq P$,
 \item[{(ii)}] $\mathrm{det}(Q)=k^*$.
 \end{itemize}

 \noindent {\bf Proof.} \rm We recall from Theorem 3.5 that in both cases $G \backslash \TTT$ is a {\it tree}. Let $\TTT_0$ be a {\it lift} of $G \backslash \TTT$ with respect to the natural projection of $\TTT$ onto $G \backslash \TTT$. Hence, by definition, $\TTT_0$ is a subtree of $\TTT$ isomorphic to $G\backslash \TTT$. It is known that there exists $e \in \edge(\TTT_0)$ for which
 $$G_e=Z.$$
 \noindent The edge $e$ naturally partitions the vertices of $\TTT_0$ into $V_1,\;V_2$, say. Let $P=\left< G_v: v \in \vert(V_1)\right>$ and
 $Q=\left< G_v: v \in \vert(V_2) \right>$. Then from standard Bass-Serre theory [Se, p.42]
 $$G=P\star_{\quad Z}Q.$$
 \noindent We can choose $\TTT_0$ so that there exists $v \in \vert(V_1)$ for which
 $$G_v=GL_2(k).$$
 \noindent In addition there exists $v \in \vert(V_2)$ for which $G_v$ is of type (a) or (c). From the descriptions of the matrices in the stabilizers of these types given in [MS4, Theorems 2.1, 2.6] it is clear that in either case
 $$\det(G_v)=k^*.$$
 \noindent (For stabilizers of type (a) we require the fact that the {\it norm} map $N_{L/k}$, where $L=\FF_{q^2}$, is {\it surjective}.)\\ \\
 \noindent For the elliptic case $(g,\delta)=(1,1)$ we can choose $\TTT_0$ so that the assertions hold. See Takahashi's paper [T], in particular [T, Theorems 3, 5].\\ \\
 \noindent For the case $(g,\delta)=(0,3)$ a detailed description of $\TTT_0$ is given in [M5, Theorem 2.26] (for the case $d=3$). The edge with trivial stabilizer is, in the notation of [M5], the one joining $\bar{\Lambda}_0$ and $\bar{\Lambda}_1[1]$. It turns out we can choose $\TTT_0$ so that the stabilizer of $\bar{\Lambda}_0$ is $GL_2(k)$ and that {\it all} other stabilizers of the vertices of $\TTT_0$ are of type (c).
\hfill $\Box$
\\ \\
\noindent {\bf (a) The elliptic case $(g,\delta)=(1,1)$.}
\\ \\
\noindent (i) Suppose that $4|q$. Then $\r(\Gamma)=\r(G)=0$.
By [MS2, Theorem 5.5] and Proposition 4.3
$$\ngncs(\Gamma)=\m (\Gamma)=p',$$
\noindent where, with two exceptions, $p'=3$. For the exceptional cases
(when $q=4$) $p'=5$.
\par
Take $M \unlhd \Gamma$ with $|\Gamma:M|=p'$ and let $N=Z.M$. Then, since
$Z.\Gamma=G$ it follows that $N \unlhd G$ and that $|G : N|=p'$. By
Proposition 4.1 and Remark 4.2 $N$ is a normal genuine non-congruence
subgroup of $G$. So $\ngncs(G)\leq p'$.
\par
Now let $H$ be any subgroup of $G$ with $|G:H|<p'$. Then
$|\Gamma : H\cap\Gamma|<p'$, and hence $\Gamma\leq H$.
We conclude then that, when $4|q$,
$$\ngncs(G)=p'.$$
\noindent With the two above exceptions $\ngncs(G)=\m(G)=3$. For the two
exceptional cases (when $q=4$) we have $\ngncs(G)=\gncs(G)=5$ and $\m(G)=3$.
The subgroup attaining the latter bound is (the {\it congruence} subgroup)
$\Gamma$.
\\ \\
\noindent (ii) Suppose now that $q$ is odd. From the rank formulae
([G1] or see [MS1, p.33]) it follows that here $r(\Gamma)=q$ and so
Theorem 5.1 applies here for the case $X=\Gamma$.  However $\r(G)=0$.
With the notation of Lemma 5.2 we define an epimorphism $\phi$ from $G$ to $\{\pm1\}$ by
$$
\phi(g)=\left\{\begin{array}{cll}
1 &, & g\in P\\[8pt]
(\det(g))^{\frac{q-1}{2}} &, & g\in Q\\
\end{array}\right.
$$
\noindent Hence there exists $N \trianglelefteq G$, containing $GL_2(k)$
for which $|G:N|=2$. By Proposition 4.1 therefore in this case
$$\ngncs(G)=\mathrm{m}(G)=2.$$
\noindent We summarize the results for this case.
\\ \\
\noindent {\bf Theorem 5.3.} \it
Suppose that $(g,\delta)=(1,1)$ and that $q \neq 2$.
\begin{itemize}
\item[(i)] If $q$ is odd
$$\ngncs(X)=\gncs(X)=\m(X)=2.$$
\item[(ii)] If $4|q$ then with two exceptions
$$\ngncs(X)=\gncs(X)=\m(X)=3.$$
For each of the exceptional cases $q=4$ and
$$\ngncs(X)=\gncs(X)=\m(\Gamma)=5,\;\m(G)=3.$$
\end{itemize}
\rm
\smallskip
\noindent
{\bf Remark 5.4.}
Less precise results appear to hold for the elliptic case when $q=2$.
When $q=2$ it is known [MS1, Lemma 3.1(i)] that, for {\it any} $A$,
$$\ncs(A)=2.$$
If $A$ is of elliptic type and has not more than two prime ideals of 
degree $1$ (in other words, if the underlying elliptic curve has at 
most $3$ $\FF_2$-rational points), then by Takahashi's results 
[T, Theorem 5] $GL_2(A)$ has a free factor $\ZZ/3\ZZ$. Consequently,
$GL_2(A)$ has a normal subgroup of index $3$, which by 
Proposition 4.3 must be genuine. So in this case
$$2 \leq\gncs(G)=\ngncs(G)\leq 3.$$
\noindent We decide the most interesting case [Se, 2.4.4, p.115].
\\ \\
{\bf Example 5.5.} Let 
$$A=\FF_2[x,y]\ \ \ \hbox{\rm with}\ \ \ y^2 +y=x^3 +x+1.$$
Then $A$ has no prime ideals of degree $1$. So by Lemma 2.5 every
subgroup of index $2$ is a genuine non-congruence subgroup.
\\ \\
{\bf (b) The case $(g,\delta)=(0,3),\; q \neq 2$.}
\\ \\
\noindent \rm (i) Suppose that $4|q$. It follows from
[MS2, Theorem 4.7, Theorem 6.2] and Proposition 4.3 that
$$\ngncs(\Gamma)=\mathrm{m}(\Gamma)= p'$$
\noindent where $p'$ is the smallest prime dividing $q-1$.
Taking $M \trianglelefteq \Gamma$ with $|\Gamma:M|=p'$ and considering
the subgroup $Z.M$ it can be shown as in (a)(i) above that
$$\ngncs(G)=\mathrm{m}(G)=p'.$$
\\
\noindent (ii) Suppose that $q$ is odd. Then as in the elliptic case from Lemma 5.2 it follows that
$$\ngncs(G)=\mathrm{m}(G)=2.$$

\noindent We summarize the results for this case.
\\ \\
{\bf Theorem 5.6.} \it
Suppose that $(g,\delta)=(0,3)$ and that $q \neq 2$. Then
$$\ngncs(X)=\gncs(X)=\mathrm{m}(X)=p',$$
\noindent where $p'$ is the smallest prime dividing $q-1$.
\rm
\\

\subsection*{6. The case $A=\FF_q[t]$}

\noindent Finally we look at the most important case $(g,\delta)=(0,1)$,
that is, $A=\FF_q[t]$. For most aspects of Drinfeld modular curves this
is the case that is by far the best understood. Ironically, for this
case we don't know $\ngncs(G)$ exactly and can only give lower bounds.
\\ \\
{\bf Theorem 6.1.} \it
Let $A=\FF_q[t]$ with $q>3$ and let $N$ be a normal genuine non-congruence
subgroup of $\Gamma$. Then $|\Gamma:N|$ is divisible by
$q|PSL_2(\FF_q)|$. In particular
$$\ngncs(\Gamma)>\nncs(\Gamma)=|PSL_2(\FF_q)|.$$
\rm
\\
{\bf Proof.}
Actually, the index of any proper normal subgroup of $\Gamma$ is divisible
by $|PSL_2(\FF_q)|$. Namely, by [MS3, Lemma 3.2] we either have $\ql(N)=A$
(and hence $N=\Gamma$) or $N\cap SL_2(\FF_q)\leq\{\pm I_2 \}$. In the latter
case $SL_2(\FF_q)/(N\cap SL_2(\FF_q))$ is a subgroup of $\Gamma/N$.
\par
So we only have to show that $q^2$ divides $|\Gamma:N|$. Assume not.
Since $A/\ql(N)$ is a subgroup of $\Gamma/N$, this implies that $\ql(N)$
has codimension $1$ (or $0$) in $A$. From the previous paragraph we already
know that $T(1)\notin N$, that is, $1\notin \ql(N)$. Hence there exists
a non-standard automorphism $\varPhi$ of $\Gamma$ with $\ql(\varPhi(N))=tA$.
But $\varDelta(t)=\Gamma(t)$ ([M1, Corollary 1.4]); so $\varPhi(N)$ is
a congruence subgroup.
\hfill $\Box$
\\ \\
The corresponding result for $q\leq 3$ requires a little bit of preparation.
\\ \\
{\bf Lemma 6.2.} \it
Let $q\leq 3$ and $A=\FF_q[t]$. Denote the commutator group of $\Gamma$
by $\Gamma'$.
\begin{itemize}
\item[(i)] $\Gamma'$ is the normal subgroup of $\Gamma$ generated by the
unique subgroup of order $q^2 -1$ of $SL_2(\FF_q)$.
\item[(ii)] If $N$ is a normal subgroup of $\Gamma$ containing $\Gamma'$,
then $\Gamma/N$ is naturally isomorphic to $A/\ql(N)$.
\item[(iii)] If $N$ is a normal subgroup of $\Gamma$ with $1\in\ql(N)$,
then $N$ contains $\Gamma'$.
\end{itemize}
\rm

\noindent{\bf Proof.}
For $q\leq 3$ the commutator of $SL_2(\FF_q)$ is its unique subgroup
$P$ of order $q^2 -1$. So $\Gamma'$ contains the normal hull of $P$.
For the converse we use Nagao's Theorem
$$\Gamma=SL_2(\FF_q)\star_{SB(\FF_q)} SB(A),$$
where $SB(R)=B_2(R)\cap\Gamma$. See for example [Se, p.88].
From this we see that the quotient of $\Gamma$ by the
normal hull of $P$ is $T(A)$. This proves (a) and (b).
\par
Part (c) follows from the simple fact that the normal subgroup of
$SL_2(\FF_q)$ generated by $T(1)$ is $SL_2(\FF_q)$ itself.
\hfill $\Box$
\\ \\
{\bf Theorem 6.3.} \it
Let $A=\FF_q[t]$ with $q\leq 3$ and let $N$ be a normal genuine non-congruence
subgroup of $\Gamma$. Then $q^2$ divides $|\Gamma:N|$.
In particular
$$\ngncs(\Gamma)>\nncs(\Gamma)=q.$$
\rm

\noindent {\bf Proof.}
If $1\notin\ql(N)$, the proof is exactly the same as for Theorem 6.1.
\par
If $1\in\ql(N)$ and $q^2\nmid |\Gamma:N|$, then by Lemma 6.2 necessarily
$\Gamma'\leq N$ and $|\Gamma:N|=|A:\ql(N)|=q$. Then there exists an
$\FF_q$-vector space automorphism $\phi$ of $A$ (with corresponding
non-standard automorphism $\varPhi$) such that $\phi(1)=1$ and
$\phi(\ql(N))=t(t-1)A\oplus\FF_q$. So $\varPhi(N)$ has level $t(t-1)A$
(and still contains $\Gamma'$).
\par
To finish the proof we verify that
$\varPhi(N)$ is a congruence subgroup by showing
$$\varDelta(t(t-1)).\Gamma' =\Gamma(t(t-1)).\Gamma'.$$
It is well-known that
$$\Gamma/\Gamma(t(t-1))\cong SL_2(A /(t))\times SL_2(A/(t-1)).$$
\noindent Under this isomorphism
$$(\Gamma(t(t-1)).\Gamma')/\Gamma(t(t-1))\cong
SL_2(\FF_q)'\times SL_2(\FF_q)'.$$
It follows that
$$|\Gamma:\Gamma(t(t-1)).\Gamma'|=q^2=|\Gamma:\varDelta(t(t-1)).\Gamma'|,$$
which proves the claim.
\hfill $\Box$
\\ \\
{\bf Theorem 6.4.} \it
Let $N$ be a normal genuine non-congruence subgroup
of $G=GL_2(\FF_q[t])$. Then $|G:N|$ is divisible by
$$
\left\{\begin{array}{cl}
q|PGL_2(\FF_q)|=q^4 -q^2, & \hbox{\rm if}\ q>3,\\[8pt]
q^2,  & \hbox{\rm if}\ q\leq 3.\\
\end{array}\right.
$$
In particular,
$$\ngncs(G)>\nncs(G).$$
\rm
\\
{\bf Proof.}
First of all, $|G:N|$ is divisible by $|\Gamma:N\cap\Gamma|$, which is
bigger than $1$ unless $\Gamma\subseteq N$. If $q>3$, then, as explained
earlier, $|\Gamma:N\cap\Gamma|$ is divisible by $|PSL_2(\FF_q)|$. So
$|G:N|\geq\frac{q^3 -q}{2}$, and hence $|G:ZN|\geq\frac{q(q+1)}{2}>q-1$.
In particular, $ZN$ does not contain $\Gamma$. 
\par
Now $GL_2(\FF_q)\cap ZN$ is a normal subgroup of $GL_2(\FF_q)$. But from
the proof of Theorem 6.1 we know that $SL_2(\FF_q)\cap ZN=\{\pm I_2\}$.
This leaves only the possibility $GL_2(\FF_q)\cap ZN=Z$. So
$$GL_2(\FF_q)/(GL_2(\FF_q)\cap ZN)\cong PGL_2(\FF_q)$$ 
is a subgroup of $G/(ZN)$. This shows that $|G:N|$ is divisible 
by $q^3 -q$ if $q>3$.
\par
Now assume that $q^2$ does not divide $|\Gamma:N\cap\Gamma|$. By the proofs
of Theorems 6.1 and 6.3 then there exists a non-standard automorphism $\varPhi$
of $\Gamma$ that maps $N$ to a congruence subgroup. By defining $\varPhi$ to
act as identity on diagonal matrices, $\varPhi$ extends to a non-standard
automorphism $\varPhi$ of $G$. As $\varPhi(N)$ contains a congruence subgroup,
$N$ is not genuine.
\par
To prove the last claim we exhibit a normal congruence subgroup of small
index in $G$ and quasi-level $tA$. By a suitable non-standard automorphism
this group can then be mapped to a normal (non-genuine) non-congruence
subgroup of the same index.
\par
If $q>3$ we can take $Z.\Gamma(t)$, which has index $|PGL_2(\FF_q)|$.
If $q=2$, then $G=\Gamma$ anyway, and Theorem 6.3 applies.
Finally, if $q=3$ we observe that the $2$-Sylow subgroup of
$SL_2(\FF_3)$ is normal in $GL_2(\FF_3)$. Taking its inverse image
under the (in this case surjective) natural map $G\to GL_2(A/(t))$,
we obtain a normal subgroup of index $6$.
\hfill $\Box$
\\ \\
{\bf Remark 6.5.}
More precisely, Theorems 6.1, 6.3 and 6.4 show that in order for a normal
subgroup $N$ to have a chance of being genuine $\ql(N)$ must have at least
codimension $2$ in $\FF_q[t]$, and hence $X/N$ must contain a subgroup 
isomorphic to $\FF_q\oplus\FF_q$.
\\ \\
We finish with a partial result on not necessarily normal genuine
non-congruence subgroups.
\\ \\
{\bf Corollary 6.6.} \it
Let $q=p$ a {\it prime}, i.e. $A=\FF_p[t]$. Then $\gncs(X)\geq 2p$,
and hence in particular $\gncs(\Gamma)>\ncs(\Gamma)$.
\rm
\\ \\
{\bf Proof.}
Let $N$ be the core of $H$ in $X$. If $|X:H|<2p$, then $|X:N|$ divides
$(2p-1)!$ and is therefore not divisible by $p^2$. So $N$ cannot be genuine,
and consequently neither can be $H$.
\hfill $\Box$
\\ \\ \\
{\bf Acknowledgements.} The second author thanks ASARC in South Korea for
support and in particular for financially supporting his research visits
to Glasgow University.
\\

\subsection*{\hspace*{10.5em} References}
\begin{itemize}

\item[{[B]}] H.~Bass: \it Algebraic K-Theory, \rm (Benjamin, New York, 1968).

\item[{[D]}] J.~L.~Dyer: Automorphism sequences of integer unimodular groups,
\it Illinois J. Math. \bf 22 \rm (1978), 1-30.

\item[{[G1]}] E.-U.~Gekeler:  Le genre des courbes modulaires de Drinfeld,
\it C.R. Acad. Sci. Paris \bf 300 \rm (1985), 647-650.

\item[{[G2]}] E.-U.~Gekeler: \it Drinfeld Modular Curves, \rm
(Springer LNM 1231, Berlin Heidelberg New York, 1986).

\item[{[HR]}] L.~K.~Hua and I.~Reiner: Automorphisms of the unimodular group,
\it Trans. Amer. Math. Soc. \bf 71 \rm (1951), 331-348.

\item[{[JT]}] G.~A.~Jones and J.~S.~Thornton: Automorphisms and congruence
subgroups of the extended modular group, \it J. London Math. Soc. (2) \bf 34
\rm(1986), 26-40.

\item[{[L]}] A. Lubotzky: Free quotients and the congruence kernel of $SL_2$, 
\it J. Algebra \bf 77 \rm (1982), 411-418.

\item[{[M1]}] A.~W.~Mason: Anomalous normal subgroups of $SL_2(K[x])$,
\it Quart. J. Math. Oxford (2) \bf 36 \rm (1985), 345-358.

\item[{[M2]}] A.~W.~Mason: The order and level of a subgroup of $GL_2$ over
a Dedekind ring of arithmetic type, \it Proc. Roy. Soc. Edinburgh Sect. A
\bf 119 \rm (1991), 191-212.

\item[{[M3]}] A. W. Mason: Quotients of the congruence kernels of $SL_2$ over 
arithmetic Dedekind domains, \it Israel J. Math. \bf 91 \rm (1995), 77-91.

\item[{[M4]}] A.~W.~Mason: On non-congruence subgroups of the analogue of
the modular group in characteristic $p$: Rankin memorial issues.
\it Ramanujan J. \bf 7 \rm (2003), 141-144.

\item[{[M5]}] A.~W.~Mason: The generalization of Nagao's theorem
to other subrings of the rational function field,
\it Comm. Algebra \bf 31 \rm (2003), 5199-5242.

\item[{[MPSZ]}] A.~W.~Mason, A.~Premet, B.~Sury, P.~A.~Zalesskii:
The congruence kernel of an arithmetic lattice in a rank one algebraic 
group over a local field,
J. Reine Angew. Math. \bf 623 \rm(2008), 43-72.

\item[{[MS1]}] A.~W.~Mason and A.~Schweizer: The minimum index of
a non-congruence subgroup of $SL_2$ over an arithmetic domain,
\it Israel J. Math. \bf 133 \rm (2003), 29-44.

\item[{[MS2]}] A.~W.~Mason and A.~Schweizer: The minimum index
of a non-congruence subgroup of $SL_2$ over an arithmetic domain.
II: The rank zero cases, \it J. London Math. Soc. \bf 71 \rm
(2005), 53-68.

\item[{[MS3]}] A.~W.~Mason and A.~Schweizer: Non-standard automorphisms
and non-cong\-ruence subgroups of $SL_2$ over Dedekind domains contained
in function fields, \it J. Pure Appl. Algebra \bf 205 \rm (2006), 189-209.

\item[{[MS4]}] A.~W.~Mason and A.~Schweizer: The stabilizers in a Drinfeld
modular group of the vertices of its Bruhat-Tits tree: an elementary approach,
\it Int. J. Algebr. Comput. \bf 23 \rm (2013), 1653-1683.

\item[{[R]}] I.~Reiner: A new type of automorphism of the general linear
group over a ring, \it Ann. of Math. \bf 66 \rm (1957), 461-466.

\item[{[Se]}] J.-P.~Serre, \it  Trees,  \rm (Springer, Berlin,
Heidelberg, New York, 1980).

\item[{[SV]}] J.~Smillie and K.~Vogtmann: Automorphisms of $SL_2$ over
imaginary quadratic integers, \it Proc. Amer. Math. Soc. (2) \bf 112
\rm (1991), 691-699.

\item[{[St]}] H.~Stichtenoth: \it Algebraic Function Fields and Codes
(Second Edition), \rm (Springer GTM 254, Berlin Heidelberg, 2009).

\item[{[T]}] S.~Takahashi: The fundamental domain of the tree
of $GL(2)$ over the function field of an elliptic curve,
\it Duke Math. J.  \bf 72  \rm(1993), 85-97.

\end{itemize}

\end{document}